# On Free Fall in the Three Body Problem

*Dedicated to Prof. R. Broucke*

by F.L.M.G Janssens[1]

ABSTRACT

The free fall of three particles under Newtonian attraction allows to illustrate some of the complexities of the general three body problem. The total collapse or singularity that occurs when starting from one of the five central configurations (two triangular and three collinear) generates periodic solutions and the singularity mimics an elastic bounce. Periodic solutions without collisions where found by Standish : three particles fall from an initial triangle to each other and without colliding, come later to rest on another triangle where the motion reverses. Singularities where the motion ends, are illustrated by equal particles starting from an isosceles triangle.

The lack of continuity in neighbouring solutions is illustrated by particles starting from a nearly equatorial triangle. Although the total energy is negative, an elliptic-hyperbolic break up of the system where all three particles go to infinity. is possible. Two particles are tightly bound in elliptic motion, their CoM recedes to infinity while in an hyperbolic motion with the third particle . The famous historic case of the Pythagorean triangle shows that such a break up may happen after a long time and several close passages. The break-up in a elliptic hyperbolic system occurs in a very short time period around a very close passage. Progress in the understanding the interactions between the particles when they are very close,can lead to sharper escape criteria.

This review suggests that for the free fall , there are only three types of final trajectories :
a) periodic with or without collisions, b) ending in a ternary collision and c) a break-up in an elliptic-hyperbolic system.

## INTRODUCTION

### Central Configurations - Triple Collisions (TC)

The initial configurations that lead to a total collapse are known as "central configurations " . In the 2BP, any initial positions of the two particles is a central configuration. The central configurations fill the configuration space. In the 3BP, the two equilateral and three collinear configurations are the only central configurations and they cover only part of the configuration space.

### A. Three points on an equilateral triangle (Lagrange solution).

When 3 masses $m_i$ ($\Sigma m_i = M$) are initially located on an equilateral triangle, their acceleration points the CoM. Indeed , the equations of motion , with the notations of [Steves and & 1977] are :

$$m_1 \ddot{\bar{r}}_1 = Gm_1(m_2\bar{\rho}_{12} + m_3\bar{\rho}_{13})$$
$$m_2 \ddot{\bar{r}}_2 = Gm_2(m_1\bar{\rho}_{21} + m_3\bar{\rho}_{23})$$
$$m_3 \ddot{\bar{r}}_3 = Gm_3(m_1\bar{\rho}_{31} + m_2\bar{\rho}_{32})$$

(1)

$$\text{where} \quad \bar{\rho}_{ij} = \frac{\bar{r}_j - \bar{r}_i}{r_{ij}^3} \quad \text{and} \quad \sum m_i = M$$

When the CoM is at rest:

$$\sum_{i=1}^{3} m_i \bar{r}_i = \sum_{i=1}^{3} m_i \dot{\bar{r}}_i = \sum_{i=1}^{3} m_i \ddot{\bar{r}}_i = \mathbf{0}$$

(2)

For particles located on an equilateral triangle : $r_{12} = r_{23} = r_{31} = r$ . Using the definition of the CoM , Eq. (1) can be written as :

$$\ddot{\bar{r}}_i = -\frac{\mu}{r^3}\bar{r}_i \quad \text{where} \quad \mu = GM, \text{ and } i = 1,3$$

(3a)

---
[1] f.anssens@ziggo.nl , retired from ESTEC/ESA

which shows that the accelerations are central and the three particles start to move in to the origin. The denominator r is not the distance to origin ($r_i$) but the side of the triangle.
By subtracting the Eq.(1) from each other, one obtains the equation of motion for the side of a triangle $\bar{r}$ :

$$\ddot{\bar{r}} = -\frac{\mu}{r^3}\bar{r} \quad \text{where} \quad \mu = GM \tag{3b}$$

which is identical to the equations for the Kepler problem (fixed origin) or the equations for the relative motion of the 2BP. The general relation between the sides ($r_{ij}$) of a triangle with masses $m_i$ at the vertices and the distance to the CoM is :

$$\begin{bmatrix} r_{23}^2 \\ r_{31}^2 \\ r_{12}^2 \end{bmatrix} = \begin{bmatrix} -\frac{m_1^2}{m_2 m_3} & \frac{m_2+m_3}{m_3} & \frac{m_2+m_3}{m_2} \\ \frac{m_3+m_1}{m_3} & -\frac{m_2^2}{m_3 m_1} & \frac{m_3+m_1}{m_1} \\ \frac{m_1+m_2}{m_2} & \frac{m_1+m_2}{m_1} & -\frac{m_3^2}{m_1 m_2} \end{bmatrix} \begin{bmatrix} r_1^2 \\ r_2^2 \\ r_3^2 \end{bmatrix} \quad \text{and} \quad \begin{bmatrix} r_1^2 \\ r_2^2 \\ r_3^2 \end{bmatrix} = \frac{1}{M^2}\begin{bmatrix} -m_2 m_3 & (m_3+m_2)m_3 & (m_2+m_3)m_2 \\ (m_3+m_1)m_3 & -m_3 m_1 & (m_3+m_1)m_1 \\ (m_1+m_2)m_2 & (m_1+m_2)m_1 & -m_1 m_2 \end{bmatrix}\begin{bmatrix} r_{23}^2 \\ r_{31}^2 \\ r_{12}^2 \end{bmatrix}$$

(4a,b)

For an equilateral triangle $r_{ij} = r$, Eq.4 simplifies to :

$$r_i^2 = \frac{m_j^2 + m_k^2 + m_j m_k}{M^2} r^2 \quad \text{with } j,k \neq i \tag{5}$$

while for equal masses the matrix reduces to : $\begin{bmatrix} -1 & 2 & 2 \\ 2 & -1 & 2 \\ 2 & 2 & -1 \end{bmatrix}$ and is proportional to its inverse.

So, substituting Eq.5 in Eq.3 :

$$\ddot{\bar{r}}_i = -\frac{\mu_i}{r_i^3}\bar{r}_i \quad \text{where} \quad \mu_i = GM_i \quad M_i = \frac{(m_j^2 + m_k^2 + m_j m_k)^{3/2}}{M^2} \tag{6}$$

which is the same as Eq.(3b) with a different µ value. Starting with zero velocity, every particle moves as if attracted by as mass $M_i$ placed at the origin. Some examples are given in the Fig.1a-c below. For equal masses, $M_i = M$ ,the total mass. The trajectories are simply the straight lines to the CoM. The triangle itself shrinks , according to Eq.(3b) while the sides remain parallel to their initial direction. The relative motion of each particle w.r.t. the other ones is along a side of the equilateral triangle which has a fixed direction. When the 3 particles have different masses, they move with different velocities on lines that meet at the Com and not at the centre of the equilateral triangle.

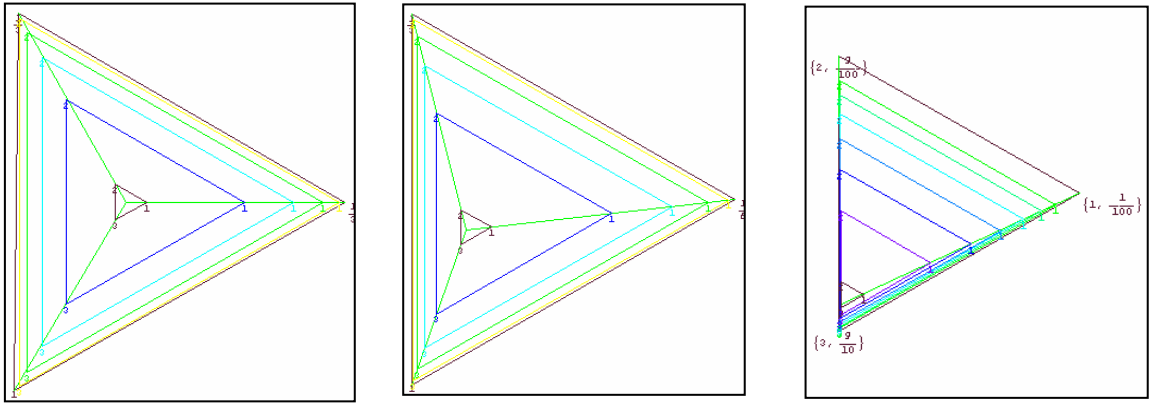

**Fig.1a - equal masses 1/3,1/3,1/3    Fig.1b - masses in ration  1/6 1/3 1/2    Fig.1c  masses 1/100 9/100 9/10**

A more rigorous treatment that the particles stay on an equilateral triangle, is given in [Roy, p.123].

When the initial velocities are not zero but compatible with a single angular velocity (any value) perpendicular to the plane of the three particles, the triangle rotates while shrinking. This is Lagrange's triangular case. Some of these trajectories are illustrated in the figures below.

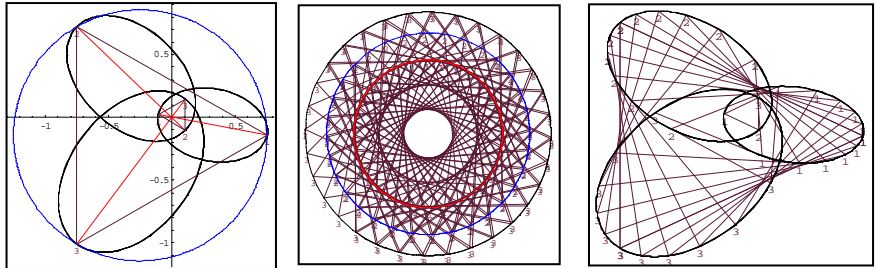

**Fig2 a-c. General Equilateral solutions of the 3BP**

The solution to Eq.(6) and (3b) is well known as it describes also the relative motion in the two body problem (2BP). Here we deal not with the elliptic, parabolic or hyperbolic solutions but with the straight line case and zero initial velocity. For one-dimensional motion, the separation x between two particles falling to each other from an initial separation $2a = r_0$ is given by :

$$\ddot{x} = -\frac{\mu}{x^3}x \qquad \mu = G(m_1 + m_2) \qquad (7)$$

The solution in terms of the eccentric anomaly is [Roy, p.97], counting the time from the apogee :

$$r = a(1 + \cos u) \qquad v = \dot{r} = -\sqrt{\frac{\mu}{a}} \tan\frac{u}{2} \qquad (8a,b)$$

$$t = n(u + \sin u) \qquad n = \sqrt{\frac{a^3}{\mu}} \qquad dim[n] = T$$

This solution is the limit of an ellipse with eccentricity 1, semi-minor axis (b) and parameter (p) zero[2]. The foci and the endpoints of the semi-major axis coincide. The velocity at the origin is infinite. The eccentric anomaly u is, in this context, called a <u>regularisation variable</u>  $du = a/n \; dt/r = \sqrt{(\mu/a)}\;/r\;dt$

When introducing u as independent variable in (7), the solution (8a,b) is, of course, recovered but not as expected from the simple linear harmonic oscillator equation, (Appendix C)

The time $T_c$ when the particles collide or arrive simultaneously at the origin is given by $u = \pi$

$$T_c = \frac{\pi}{\sqrt{2}\;2}\sqrt{\frac{r_0^3}{\mu}} \quad or \quad T_c = \pi\sqrt{\frac{a^3}{\mu}} \qquad a = r_0/2 \qquad (9)$$

The motion of each of the 2 particles to their fixed CoM is, of course, a simple scaling of Eq.8,9 by using the appropriate value of $\mu$. (Appendix A). Eq. 7 is also applicable for the side of a triangle and for each particle in the 3BP, with appropriate $\mu$ values : $\mu = GM$ for the side of a triangle and $\mu_i$ as given by Eq.6 for a single particle.

At the collision instant $u = \pi$. In the neighborhood $\pi + \varepsilon$ we have :

$$r = a(1 - \cos\varepsilon) \approx a\frac{\varepsilon^2}{2} \qquad \Delta t = n(\varepsilon - \sin\varepsilon) \approx n\frac{\varepsilon^3}{6}$$

and eliminating $\varepsilon$ : $\qquad \boxed{r \approx \sqrt[3]{\frac{9}{2}} \sqrt[3]{\mu} \; \Delta t^{2/3}} \qquad (10)$

which is the standard result for the behavior of r in the vicinity of this singularity. The behavior of r is independent of the initial distance and depends only on $\mu$. For a Newtonian law of attraction, $\mu$ has dimensions $L^3 T^{-2}$, so the only possibility to obtain a length (r), is to have an exponent 2/3 which follows also from the calculations. For the velocity, we obtain in a similar way :

$$\boxed{v \approx \frac{\sqrt[3]{\mu}}{\Delta t^{1/3}}} \qquad (11)$$

---

[2] Such solution is only a conic section when degenerate planar cones are included. The solution can not be written in the true anomaly which takes only the values $0, \pi$.

The expansion of r contains $t^2$ so it is invariant under a time reversal. The series expansion of v depends on $t^1$ and changes sign. If the global periodic solution were unknown, these properties show that the singularity at the origin can be described as an elastic bounce and that the motion can be continued as well in the 2BP as in the 3BP. A free fall from a central equilateral configuration leads to a periodic motion.

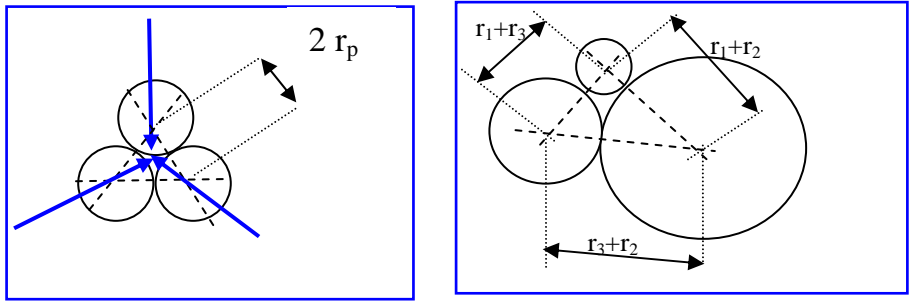

**Fig. 3a-c. Ternary collision, equal masses.; unequal masses .**

The total collapse (TC) can be considered as three elastic collisions between the particles taken two by two. It remains a fact that conservation of linear momentum and kinetic energy are not sufficient to determine the motion after a triple collision but the variable u regularizes the motion of <u>all three particles.</u> A TC does not necessarily imply that the motion ends.

Notice also that in the <u>vicinity of the collision, the quantity $rv^2$ is constant </u>and equals =2 $\mu$ from Eq.(8) for r and v. From Eq.(10,11) the constant value is $(\sqrt[3]{9/2})\,\mu = 1.65096\,\mu$. The constant 2 from Eq.(8) is the correct value as the cancellations that occur in this calculation are lost in the derivation from Eq.(10,11)

For arbitrary point masses, the CoM can be any point inside the initial triangle, The whole surface of the initial triangle, sides excluded, covers all possible locations of the CoM. When the occurrence of the triple collision is formulated in the relative motion, the side of the triangle has to become zero. As the $\mu$–value for a side of the triangle depends only on the total mass, so does the time $T_c$ from rest to the triple collision. In the cases of Fig.1, M=1 and G = 1 : $T_c$ = 2.531895753. This time is valid for any value for any non-zero value of the point masses which sum up to 1.

      The fact that a triple collision occurs for any value of the mass ratio's is far from trivial when we give each particle a radius $r_i$, related to its mass, for instance $m_i = 4/3\pi\, r_i^3$. A simultaneous collision can only occur when the distance between each of the particles is $r_i + r_j$ (Fig. 3b). They make a triangle with sides {$r_1+r_2$, $r_1+r_3$, $r_3+r_2$,}. A theorem of Sundman (1912) states that a non-collinear triple collision is only possible when the particles approach each other asymptotically on an equilateral triangle. This reuires that the radii $r_i$ are equal. Hence, in a finite picture, unequal masses must be represented by different densities and an equal radius. Notice that the particles do an head-on collision and *not an eccentric collision*, as the sides of the triangle remain parallel to themselves. Eccentric collisions can always be made to disappear by reducing the radius.

**B. Collinear configuration**

      This problem was adressed by Euler (1765) to illustrate the difficulty of the general 3BP. He derived the quintic equation that must be satisfied to have proportional solutions for the separations x,y between the points (relative motion). Roy [pp.121-125] gives a unified treatment of the equilateral and the rectilinear case. In both cases the figure made by the 3 points remains similar to itself during the motion and as a consequence each point separately conserves its total energy . In the rectilinear case, the sides of the degenerate triangle {x,y,x+y} can not be equal but they remain proportional. The distances to the CoM {$\xi_1,\xi_2,\xi_3$} remain then also proportional, of course with a different factor. In the free fall case these solutions have a triple collision at their CoM.

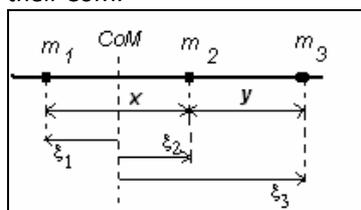

**Fig 4.- Rectilinear Case**

As already shown by Euler, the relative motion is also described by an equation of the same type as Eq.(7)
For the relative motion, $\mu_{eq}$ is :

$$\mu_{eq} = GM\left(2 + \frac{1}{(n+1)^2} - \frac{1}{n^2}\right) \qquad (12)$$

where n= y/x = $(\xi_3 - \xi_2) / (\xi_2 - \xi_1)$ is the only real root of the quintic equation:

$$(m_1 + m_2)n^5 + (3m_1 + 2m_2)n^4 + (3m_1 + m_2)n^3 = (m_2 + 3m_3)n^2 + (2m_2 + 3m_3)n + m_2 + m_3 \quad (14)$$

In the collinear case, $\mu_{eq}$ is smaller as in the equilateral case and we have no explicit expression as $\mu_{eq}$ contains the root n. The motion of each particle to the CoM is also given by a similar equation with an appropriate $\mu$ value. The $\mu_i$ values for each of the $\xi_i$ follow from the relations :

$$\xi_1 = -x\frac{m_2 + m_3(1+n)}{M} \qquad \xi_2 = x\frac{m_1 - m_3\,n}{M} \qquad \xi_3 = x\frac{n\,m_2 + m_1(1+n)}{M} \qquad (13)$$

Fig.5a-c show an example of a collinear free fall. So also in this case, the triple collision is an elastic bounce and the motion is periodic as already shown by Euler. In a finite picture of the total collapse, (particles with a radius $r_i$), only the middle particle is hit by the 2 outer ones, not all three particles touch each other.

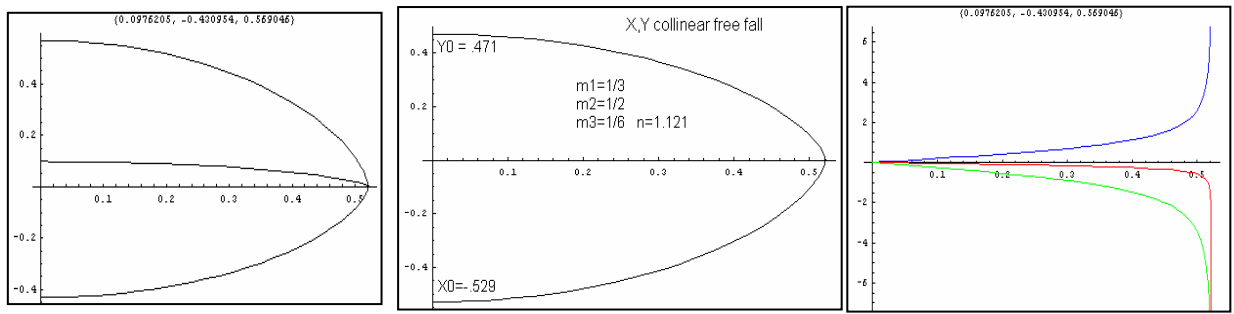

**Fig 5 a-c Rectilinear free fall : $m_1$= 1/3  $m_2$=1/2  $m_3$=1/6. Time to triple collision .5213122    n = y/x = 1.1212**

For $m_1$=$m_3$ , (n-1) factors in eqn.(14) and n=1, is the only real root. In that case or x= y , and the configuration is symmetric. Any mass $m_2$ is unmovable at the CoM. Eq.12 shows that $\mu_{eq}$ = 5/4 GM. So also in this case, the time to the TC depends only on the total mass as in the equatorial case. When $m_1 \neq m_3$, the solution n = y/x depends also on the mass ratios and not only on the total mass.
Three equal masses is a special case of $m_1$= $m_3$ .

## C. Summary  Central configurations

The triple collisions (TC) of the central configurations can be described by elastic bounces exactly as in the 2BP problem and the corresponding motion is periodic. In both cases, (collinear and equilateral), the assumption that the separation between the particles keeps a constant ratio during the motion makes the simultaneous regularisation and even an analytical solution possible. A constant energy can be ascribed to each of the particles separately. When three particles start falling from an equilateral triangle, the time to total collapse depends only on the total mass. For a collinear configuration, this property holds only when the outer masses are equal. For an arbitrary (non-degenerate) triangle, it is not possible to adjust the masses such that the particles start falling to the CoM . This emphasizes the special nature of the central configurations. The fact that the motion can be continued after a TC , is probably only possible for the central configurations.
A free fall has zero angular momentum in both the 2BP and 3BP. It turns out that a zero angular momentum is a necessary condition for a total collapse[3]. For general initial positions, the possible trajectories can be almost anything in the 3BP as will be illustrated in the sequel.

---

[3] Theorem Sundman

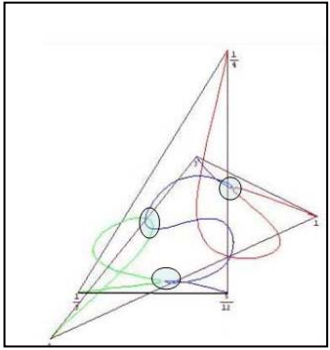

**Periodic Solution Standish (no collisions)**

Standish [1970] [4] found initial conditions for the mass ratio's { 3 : 4: 5} and a triangle close the phytagorean triangle with sides 3-4-5 considered in Burreau's problem (See further) . The initial triangle has not sides in the proportion { 3 : 4: 5 } but { 3 : 4.689 : 5.614 } and the initial right angle changed to 91.061°. These initial conditions give a periodic solution *without collisions.* The 2 masses {3:5} on the base have 2 close encounters. After a third close encounter between the masses {4:5} all three particles become simultaneously at rest on a different triangle. The motion is reversed and this time is half a period. Fig 22 shows the two triangles between which the particles oscillate.

**Fig.22 - Standish triangle.**

**Discussion Singularities**

The TC's occurring in the central configurations are probably the only ones were the motion can be continued. It was shown by [Siegel][5] , that it is *possible* that the motion cannot be continued past a TC which is then an essential singularity. Whether this situation occurs or not is not always clear for a particular TC in the 3BP[6].

When the particles remain at a TC, the potential energy is ∞ and the Kinetic energy zero. So there is a paradox as the energy must be conserved for all time. At the instant a TC occurs, the finite total energy E takes the form :

$$E = T + V = +\infty - \infty = E_0 \qquad (14)$$

*as the kinetic energy T > 0 and the potential energy U < 0. This indeterminacy **must** resolve to the initial finite value (negative and only potential for free fall) as the system is conservative.* The 3 particles can then considered as 1 particle with the total mass, and some form of internal energy.

Eq.(14) that both the Kinetic and potential energy go to infinity *at the same rate* . Let t* the instant of a total collapse then

$$\lim_{t \to t^*} \frac{T}{|V|} = 1$$

In fact, this is exactly the meaning of the relation r v² = 2 μ mentioned above as an equivalent form is

$$\frac{1}{2} v^2 = \frac{\mu}{r}$$

Note that in the 2BP, the ratio T/|U| is 1/2 for a circular orbit and varies between (1+e)/2 and (1-e)/2 for an elliptic orbit where e = eccentricity.

A general expression for the behaviour in the vicinity of a TC[7] is :

$$\mathbf{r}_i(t) = t^{2/3} \sum_{k=0}^{\infty} \mathbf{a}_{ik} t^{\alpha(m_i) k}$$

"*where α(m$_i$) is a non-constant algebraic function of the masses. When α is irrational, the motion cannot be continued past the TC*" . It seems to be the typical case, that **r**$_i$(t) can not be represented as ~ t$^\nu$ for some exponent ν, it is not necessarily true for all cases.

Next we consider particles with equal mass falling from any isosceles triangle to illustrate the occurrence of other types of TC's. Notice the this case includes the central configurations discussed above. The theoretical work of Simo and Susin[8] summarizes the known results on the possible continuation or ending of the motion.

---

[4] E.M. Standish, "New periodic orbits in the general problem of three bodies",Proceedings of Symposium Sao Paulo Sep 69, Reidel Publishing 1970..
[5] Der Dreierstoss, Ann.Math. 42,127-168, 1941
[6] Simo has obtained results for other potentials (logarithmic) and the many body problem.
[7] Arnold, Kozlov, Neistadt ," Mathematical aspects of Classical and Celestial Mechanics", Springer 2006

## Three points with equal masses on isosceles triangles

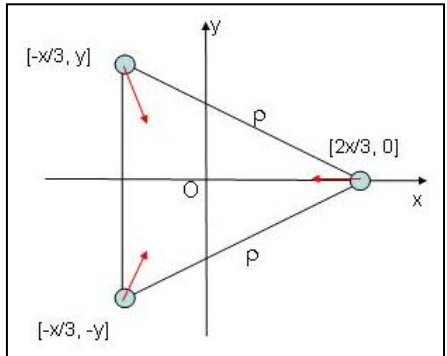

**Fig 6. Isoceles case equal masses**

The three particles are placed on the vertices of an isosceles triangle. When the two masses on the base are equal, their motion is symmetric w.r.t. to the bisector from the top on which the third particle does a one-dimensional motion. We discuss only the case where the 3 particles have equal masse m. It is clear that in this case, there are always many binary collisions. We will see that TC's occur for some initial triangles.
As shown in Fig.6, we take an axis system with the origin at the CoM. The top (vertex) has coordinates [2x/3] , the base points [-x/3, ± y]. With these notations, the immobility of the CoM is taken into account[9].The particle falling from the top moves on the x-axis.

### EoM isosceles case

The starting configuration is fully defined by the coordinates x,y, (Fig.6) and we have two equations of motion:

$$\ddot{x} = -3\frac{\mu}{\rho^3} x \quad \text{with} \quad \rho^2 = x^2 + y^2, \mu = Gm,$$

$$\ddot{y} = -\frac{\mu}{4y^2} - \frac{\mu}{\rho^3} y$$  (15a-c)

$$E = m\left\{ \frac{\dot{x}^2}{3} + \dot{y}^2 - \frac{\mu}{2|y|} - 2\frac{\mu}{\rho} \right\}$$

Eq.15b differs from Eq.7 that describes the rectilinear binary collision, by the second term that depends on x via ρ. So it is not obvious that these binary collisions have the same properties as the collisions we met in the central configurations. However, when y →0, the dominant terms of Eq.15 are :

$$\ddot{x} = -3\frac{\mu}{x^2} x \quad \text{and} \quad \ddot{y} = -\frac{\mu}{2y^2}$$  (16)

which have the same structure as Eq.7. From the y equation follows that at the limit, y satisfies Eq.10.
A change of independent variable as in Eq.8 , **du = Cst. / r(t) dt** , can not regularize both equations 15a-b simultaneously. The system contains also 2 types of singularities :

**BC => { y=0 ,x≠ 0}** $\quad E_0 = E_{BC} = \{+\infty - \infty\}_{y=0} + m\left( \frac{\dot{x}^2}{3} - 2\frac{\mu}{x} \right)$  (17)

**TC => {y=0 and x=0} equivalent with ρ=0** $\quad E_0 = E_{TC} = \{+\infty - \infty\}_{y,x=0}$  (18)

In the BC's, the energy splits up in a part for the 2 base particles and a part for the top particle.

### Discussion of numerical experiments

Tanikawa [10] has investigated free fall trajectories for equal masses numerically in a series of papers.
We discuss some of these results for an isosceles triangle according to the value of the vertex angle α :
0 < α ≤ 180°,
α = 0 is excluded as then the 3BP degenerates to a 2BP with masses in the ratio { 1:2} .
α = 60° is the equilateral central configuration
α = 180° is the collinear central configuration

---

[8] Simo, "Analysis of triple collision in the isosceles problem", in Classical Mechanics and dynamical systems, Editors Devaney,Nitecki ,Marcel Dekker NewYork, 1981
[9] These coordinates are referred to as Jacobi coordinates in "Triple Collision in the Planar Isosceles 3BP", R.L. Devaney, Inv. Math. 60, 249-267, 1980
[10] Tanikawa, Umehara e.g..: "A search for collision orbits in the free fall Three-Body Problem I, II. CM Vol 62 , 76. Theoretical investigations by Devaney and Simo.

CASE I: The vertex angle is **smaller** as 60° (Fig.7a-c).

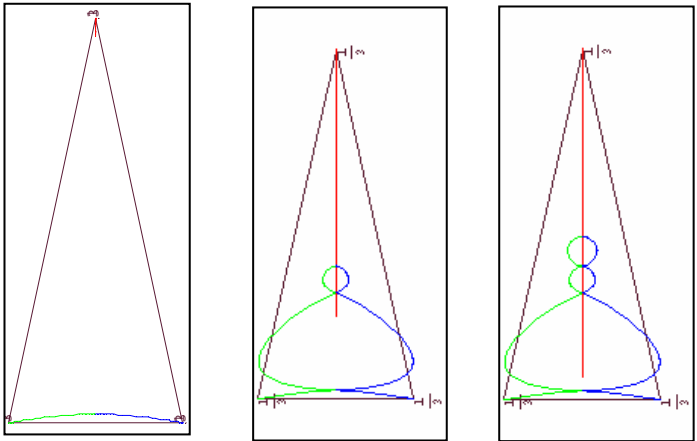

The two masses on the base make a binary collision *before* the mass moving along the symmetry axis reaches that point. The ratio of the mass at the top over two equal masses at the base cannot be adjusted to make the first binary collision a triple collision. Only when this ratio goes to infinity, the three particles collapse at the unmovable top.

**Fig.7a - c . Isosceles vertex angle 25° < 60°   (x- axis vertical)**

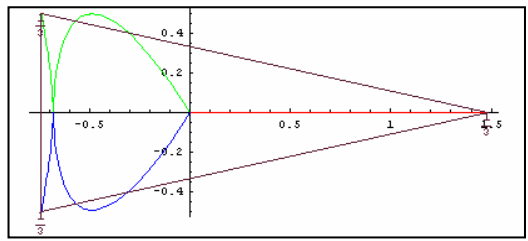

For equal masses, we can adjust the vertex angle such that a triple collision occurs on the second (third,....) binary collision of the masses starting from the base of the triangle. The first such solution occurs for $\alpha = 25.3663°$ and is shown in Fig.8. The top particle has not changed direction, and the 3 particles move meet at the CoM. The second such solution (2BC's - 1 TC) occurs for $\alpha \cong 17°$. The numerical work of Tanikawa shows an infinite number of such TC's when $\alpha \rightarrow 0$ and denotes this sequence as $T_i$ i=2,∞. More sequences of TC's are generated when reversals of the top particle are considered. In fig.8, it looks as if the TC is approached along an equilateral triangle. However, it is only at the instant of the TC when the length of the sides is zero , that the configuration is exactly equilateral. This type of TC is most likely an essential singularity.

**Fig.8 - TC for  $\alpha = 25.3663°$ , 1 BC- 1TC**

CASE II : The vertex angle is **larger** as 60°.
 the vertex particle, falling from the top, passes through zero first. At that instant we have a collinear configuration on the y-axis. The two base particles continue to move in the positive x-direction and collide eventually on the x- axis. The vertex particle continues to move in the negative x-direction while slowing down. At the instant its velocity is zero, the velocity of the particles (2,3) has to be perpendicular to the x-axis. From then onwards, all three particles move again towards each other. Contrary to case I, there is at least one reversal of direction of the vertex particle before the occurrence of a TC. Now, *we can adjust the vertex angle such that there is a TC after the **first** reversal of the mass moving on the x- axis and the first encounter of the base particles on the x-axis*. For equal masses, this angle is 140.05955° **[Tanikawa]**  [Fig 8a]. By increasing the vertex angle further, no other solutions where found in **[ Tanikawa ]** for TC 's occurs at after the second, third,....reversal of direction of the particle on the symmetry axis or after the second, third, BC of the base particles.. A priori, one expects such a series of TC's as in CASE I

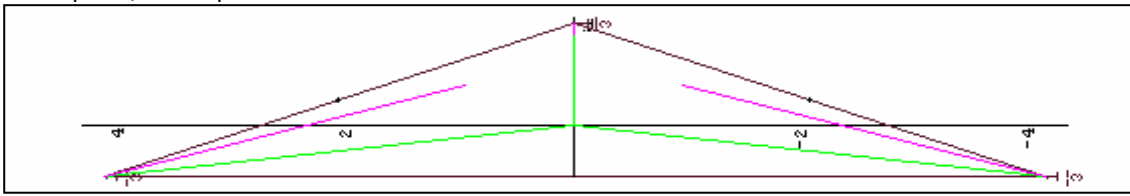

**Fig.8a Vertex angle = 140.05955°;  purple: initial acceleration;Triple collision at t=15.24696326217**

Fig. (8b-c) show TC for $\alpha$ = 140.05955° .The collinear configuration occurs at t = 12.6211 when the vertex particle is at the CoM and the two others on the y-axis.  Particle 1 comes to rest at t = 14.6194.  At that instant, the vertex angle is 130.68° and the velocity of the base particles is on the y-axis [Fig.8b]. The triple collision occurs at  t = 15.24696. Fig.8c shows how the configuration becomes equilateral when approaching the triple collision. The computations show that the configuration is only equilateral at the exact moment of the total collapse.  Fig.8c shows that with these initial conditions, the trajectory of the base particles gets again a curvature just before the triple collision.

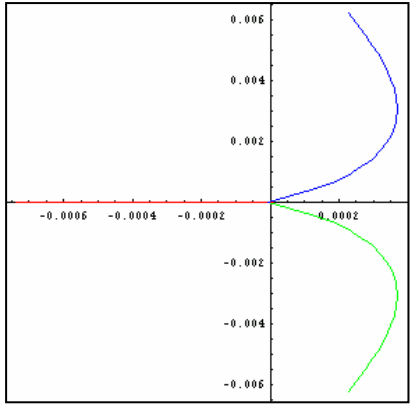
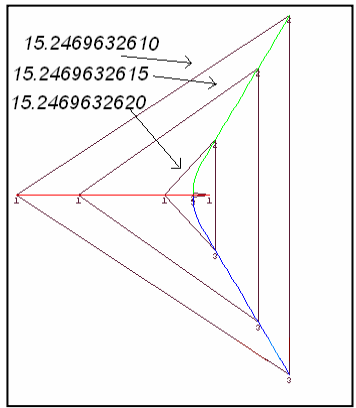

**Fig. 8b- From t=14.62 to 15.247**  **Fig.8c- Approaching the triple collision**

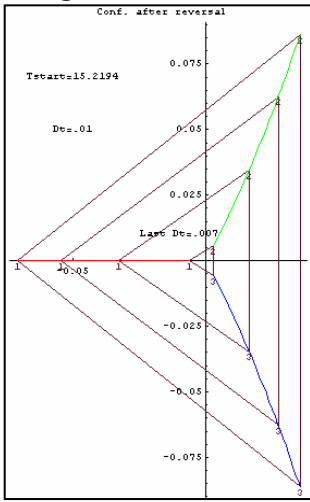
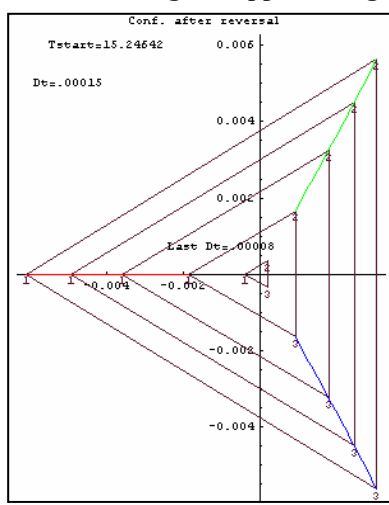
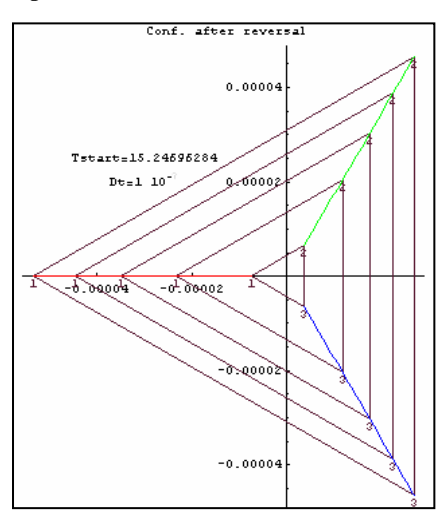

**Fig.8 d-e-f Zooms to the TC**

CASE III . Vertex angle ~180°

When the vertex angle approaches 180°, the top particle hardly moves about the CoM. When the vertex partcle is exactly is 180°, the triangle degenerates to a line segment and there is a collinear TC whjch can be continued.

Sequence of TC's

The properties of the trajectories mentioned above, suggest a systematic search for triple collisions for isosceles triangles. When the particle moving on the symmetry axis is at rest, there is a velocity along the base for the two other particles such that a triple collision occurs on the next arrival of the base particles on the symmetry axis. This velocity is outward or inward according as to the vertex angle is smaller or larger as 60° (Fig.9)

This construction defines a one-parameter family of initial conditions that leads to a TC : For a given height of an isosceles triangle with the top at rest ,there are velocities along the base line corresponding to the chosen vertex angle, such that we have a TC.

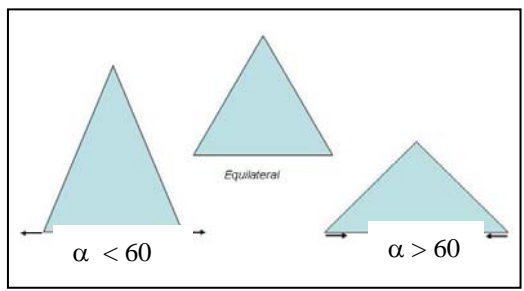
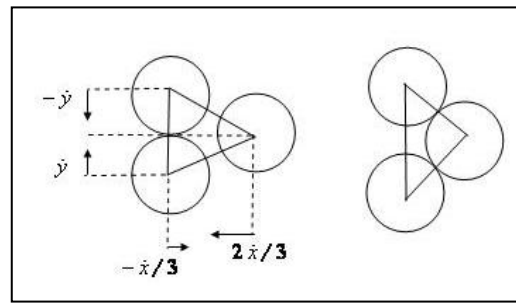

**Fig.9 - Approaching a TC when vertex at rest**    **Fig.9 a- Base Velocities for triple collision**

The case α = 140.059°) is then found by starting from an initial configuration with vertex angle 130.68°, the top particle being at rest and the base particles have an outward velocity along the y-axis of .878995.

When these isosceles configurations evolve to the TC, the configuration deforms to become an equilateral triangle at the instant of the TC. The vertex particle moves on a straight line and the two base particles follow a curved trajectory. Fig.(10a,b) show the approach to TC for α = angle 90° , base length = 3, masses 1/3 and outward velocity of the base particles $v_{y0}$ = .2204652. The ternary collision occurs at t=3.0722296.

The triangles in Fig (10a) are at $3.10^{-12}$, $2.10^{-12}$, $1.10^{-12}$, $.25\,10^{-12}$ before the collision. Fig.10b shows the calculated side length r(t) and the length given by Eq.(10). From the figures 10 above, it is tempting to conclude that that TC has the same properties as the TC's in a central configuration. However, the trajectories of the base particles do not approach the TC along a straight line of finite length. The configuration triangle continues to deform up the exact point of the TC. The velocities of the base particles is at no time directed exactly to the CoM and the relative velocity to the vertex particle is not exactly aligned with the side of the triangle.

We have a one-parameter family of initial conditions { α, $v_y(α)$ } leading to a TC that is approached in an equilateral configuration. The $t^{2/3}$ behaviour near the collision does not define the member of this family and higher derivatives must be taken into account. This is, at least, one aspect of the logarithmic singularity : the singular point can be reached from many different trajectories and the information along which trajectories the particles collided in the singular point is lost.

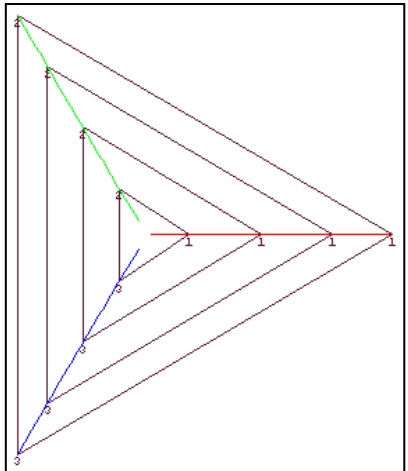 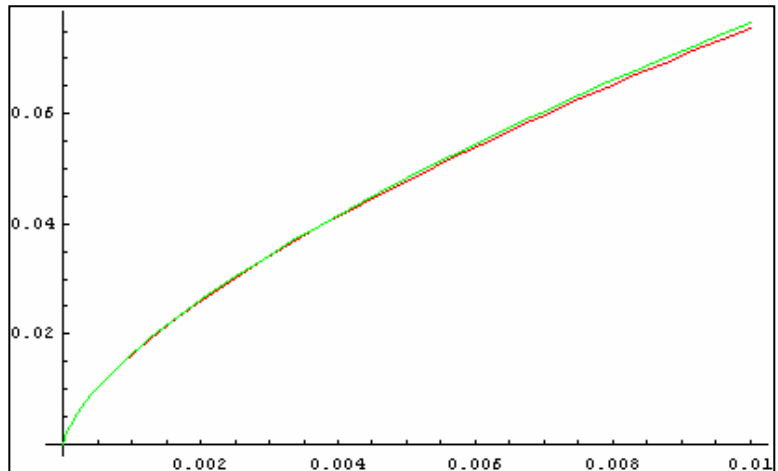

**Fig.10a - Approach TC -isoceles isoceles 90° $v_{0y}$ =.22303**   **Fig.10b- Comparison approach r(t) TC with 1. 651 $t^{2/3}$**

### Three points falling from a nearly equilateral triangle

Falling from a nearly equilateral triangle illustrates two of the astonishing properties of the general 3BP. The first property is referred to as "chaos" or chaotic behaviour and means that the slightest deviation from a set of initial conditions may result in completely different trajectories. There is, in general, no continuity in the solutions with respect to the initial conditions. For the cases in Fig. 10a-b, the lower left point (2) is shifted by .01 horizontally inwards and the upper left point (4) horizontally outwards.

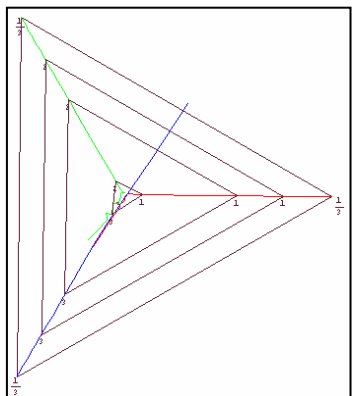 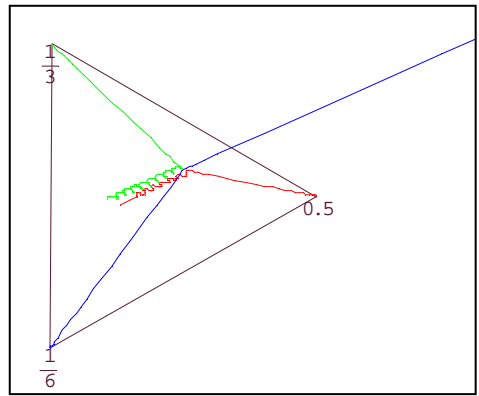

**Fig.10a - equal masses 1/3**   **Fig10b. - masses 1/6 1/3 1/2**

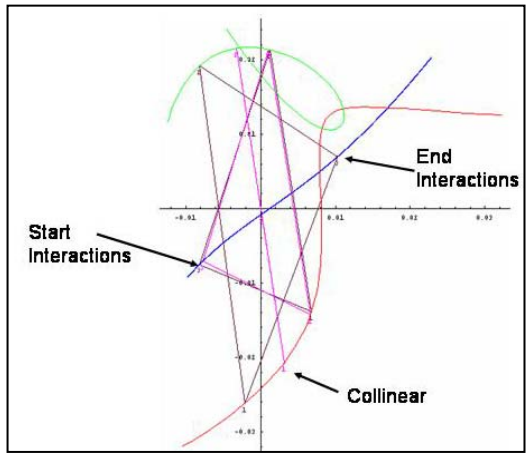 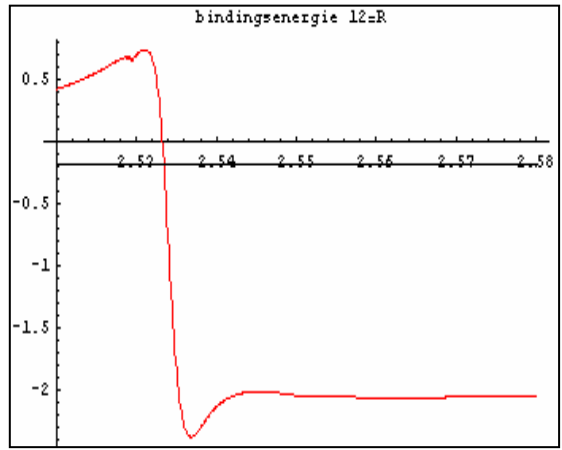

**Fig.10c - Interactions for equal masses**          **Fig 11. - Bindingsenergy 1-2**

The second, property is that with a total negative energy, the three particles can go to infinity via the split up in two asymptotically binary systems of which one is hyperbolic. This is a consequence of the interactions between the particles that take place when the triple collision is replaced by a very close encounter. In the case of Fig.10b, particle 3 passes through particle 1 and 2 which are at that instant very close to each other (collinear configuration). In the neighbourhood of this collinear configuration (Fig.10c) the system breaks up in two binary systems :  an elliptic one (particles 1,2) and an  hyperbolic system (Com of 1,2 and 3) . In these examples the split up happens in the first close encounter.

The table below shows some interesting events during the interactions of the equal mass case.  At the start of the interactions, the accelerations of 1,2 on 3 have a component increasing the velocity of mass 3. As the velocity of $m_3$ is not perpendicular to line1-2, this situation lasts until the velocity of $m_3$ makes a right angle with the line to one of the particles 1,2 (1 in this case). This happens necessarily before the collinear configuration. Simultaneously, as shown in Fig.11, the bindingsenergy of masses 1-2 decreases very fast . At t=2.53337, it becomes more negative as the total energy ($E_{tot}=E_{pot\ at\ 0}$= -.192 448). In fact, it decreases further to -2.368 at t = 2.53684. From then onwards, the bindingsenergy stabilizes on -2.06 which means that the capture is definite.

| | |
|---|---|
| 2.30841 | bindingsenergy 1-2 becomes ≥ 0 |
| 2.5294 | Min.dist 1-2 |
| 2.53250 | Quick change of bindingsenergy 1-2 starts |
| 2.53350 | where v3 is parallel to v3 at the end of the interactions |
| 2.53337 | minimum inertia |
| 2.53346 | bindingsenergy =  total (<0) energy , so here the system {1-2, 3 } becomes hyperbolic |
| 2.53371 | ang.mom. 3 = 0  (v3 // r3) |
| 2.5347 | v3max |
| **2.53494** | **collinear  configuration** |
| 2.53499 | minimum r3 toCoM, min e3 =2BP {1+2, 3}  (-199.93) |
| 2.53684 | minimum binding energy 12  (-2.368) |
| 2.53728 | the energy of the 2BP {m3, m1+m2} becomes   >  0 |
| 2.541 | ~ end of strong interactions |
| 2.55189 | NEXT MIN OF 12 (.0225 AFTER PREV./ THEN SEP. .034) |
| 2.56186 | second slow min. bindingsenergy 12 (stabilizes later on -2.060) |
| {2.585889 } | one of min dist r12 |
| {2.612995} | ""                  "" |

**Table A  - Special events during interaction**

Fig.12  shows the relative motion of the particles 1-2 about their CoM , after capture. The energy of the system CoM12, 3 becomes positive at t=2.53728. The example of Burrau (Phythagorean problem), described in the next section,  shows that such exchange of energy does  not necessarily happen in the first close encounter.

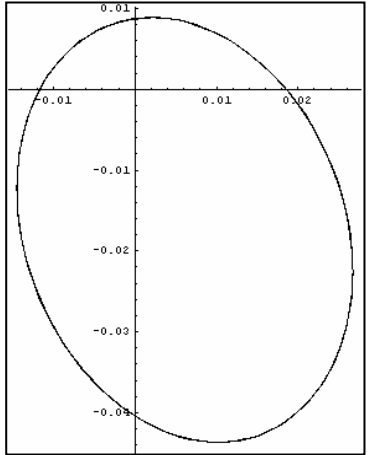 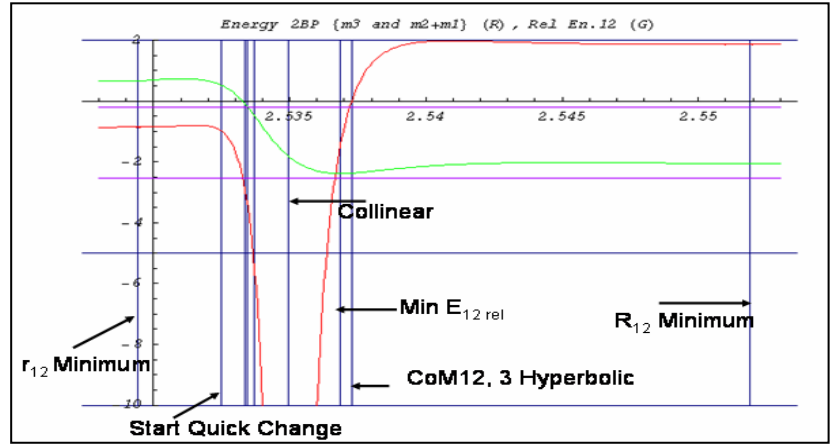

**Fig.12 - Relative motion of 1-2 after capture**   **Fig.13 - Special points Energy during interaction**

Finally , Fig.14 shows the decrease of the excess energy to the energy in the two binary systems after capture.

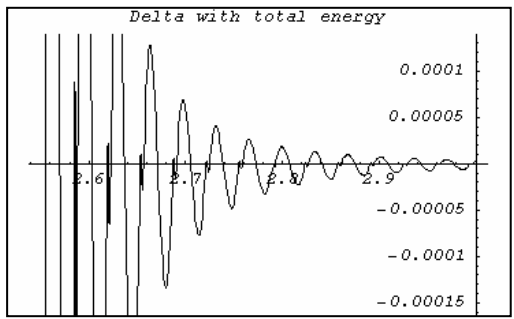 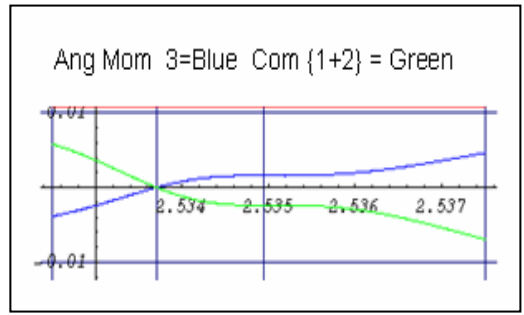

**- Fig.14 ΔExcess energy to two binary systems**   **Fig.15 - Ang. Mom. during interactions**

Fig.15 shows the angular momentum during the interactions. The trajectory of particle 3 has an inflexion point (Fig.10c). Its velocity points initially to the left of the fixed CoM of the 3 particles. When it points exactly to the CoM, its angular momentum becomes zero and then changes its sign. As noticed in Appendix B,a formal split up of the total angular momentum in two binary systems is always possible and does not imply that both these angular momenta are constant. The calculations show that both these value stabilize on ± .05 after the interactions.

This examples illustrates clearly the lack of continuity in neighbouring solutions in the 3BP (chaos). When the change in initial locations define a similar triangle or a rotated triangle, the trajectories are unaffected. This latter type of new initial conditions do not define a domain in the configuration space but only in a subspace.

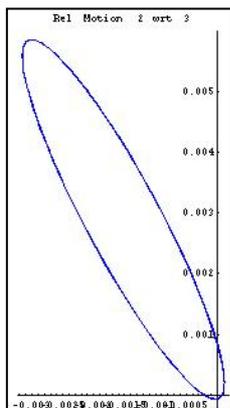 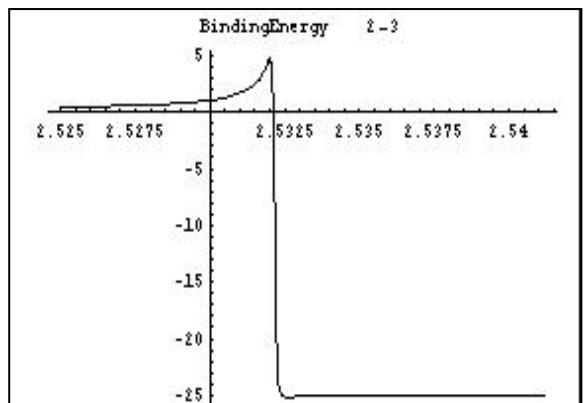

**Fig 16 a-b - Unequal masses after ejection. Binary and bindingsenergy**

A similar picture in shown in Fig 16a-b for masses {1/2, 1/3} the binary made up by and the bindings energy. After the collinear crossing, the system splits up in 2 binaries that recede to infinity. The single particle recedes faster and it looks as if that particle is ejected.

It remains challenging to find an impulsive model that summarizes the interactions. Appendix E gives starting point for the linearized equations about a BC.

Sufficient conditions for escape are e.g. given by Standish [1972][11]. Part of these conditions is always that the separation between the single particle and the CoM of the elliptic binary exceeds a value r*. In the free fall equilateral case, this distance equals the initial side of the triangle. Fig.7a-b, show that the split-up is well established immediately after the interactions of the 3 particles when they pass the global CoM. The known sufficient conditions can only confirm this particular split up when the single particle is separated from the CoM of the two other ones by at least the initial side of the triangle.

**Burrau's problem - particles on a Pythagorean triangle**

In the previous section, a single close passage caused the split up in two binary systems. The next example which has historical importance, shows that many close passages may occur before a split up takes place. In 1893[12], E.Meissel started numerical work on the Pythagorean problem: Three particles with mass 3,4,5 are on the vertices of a Pythagorean triangle opposite the sides with length 3,4,5 and start with zero initial velocity falling to each other. Meissel expected that these initial conditions would give a periodic solution. He integrated the equations numerically until just after the first near-collision ($\tau$= 2). C.Burrau[13] integrated also this case ($\tau$= 3.17) and published his results in 1913. The question if the solution is periodic was only settled in 1967 when V. Shebehely & C.F. Peters[14] continued the integration much further ($\tau$= 67) . By that time, it is clear that the three particles go to infinity, after a series of close encounters. The particles with masses 4,5 stay together in a binary system that recedes to infinity opposite to the third particle with mass 3. The conjecture of Meissel turned out to be false although at $\tau$ = 31, the 3 particles are close to their initial positions with small velocities.

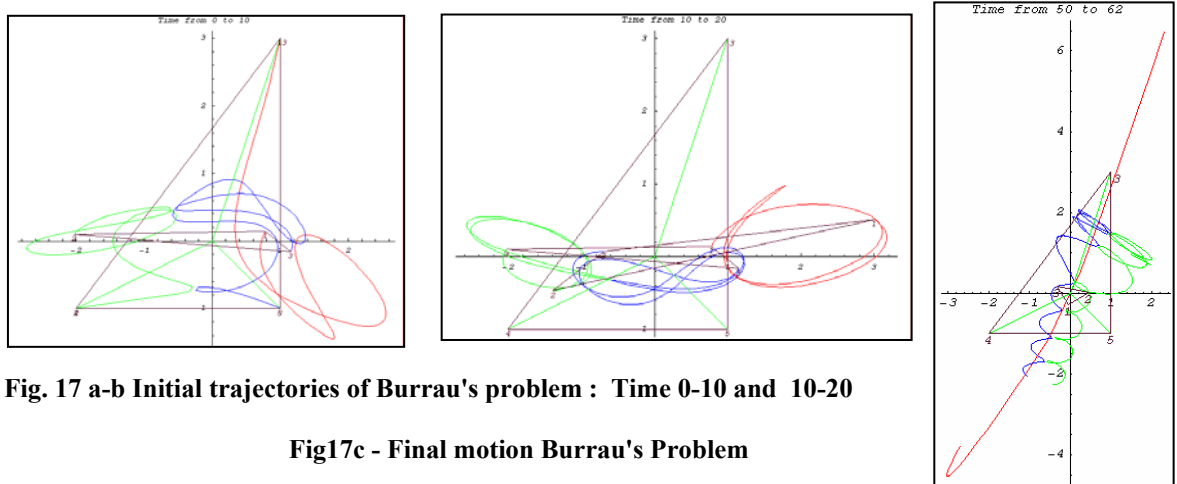

**Fig. 17 a-b Initial trajectories of Burrau's problem :  Time 0-10 and  10-20**

**Fig17c - Final motion Burrau's Problem**

The initial trajectories of the three particles is given in the figures 17a-b . The only qualifier that comes to the mind is "a hump of spagetti". The system breaks up in two binary systems around t = 59.406. ( FIg. 17c) . The particle with mass 3 has crossed the line joining the other particles 15 times before it is ejected by the system formed by the two other ones.

Fig.18 shows the bindingsenergy between the masses 4 and 5. Quick changes occur at the times of the collinear configurations, but the averaged change can go either up or down. Three crossings after t= 40 bring the bindingsenergy close to the total energy. After the crossing  at t = 46.537 , the particle with mass 3 went far out -side the initial triangle . At t=50, it is at the lower left of Fig.17c, at t = 53 it reaches a maximal distance of 5.5 from the origin and falls then back (Fig.17c,16) to pass with high velocity between the other 2 particles. The last crossing  at t = 59. 406 after which the particle is ejected. It is the first crossing where the bindings energy comes below the total energy.  Fig 19a shows that the excess energy $\Delta$ , [ Appendix B,] returns only once to 1% of the total energy at t = 59.95.

---

[11] E.M. Standish, "Sufficient Conditions for Escape in the three body problem", CM, 4(1971), 44-48

[12] Ernst Meissel and the pythagorean problem,  J.Peetre, Draft  1997

[13] Numerische berechnung eines spezialfalles des Dreikorperproblems., C.Burrau, Astron. Nachr. 195 (1913)

[14] Complete Solution of a general problem of three bodies., V.Shebehely and C.Peters, The Astonomical Journal, Vol.72,No7,(1967),pp.876-883

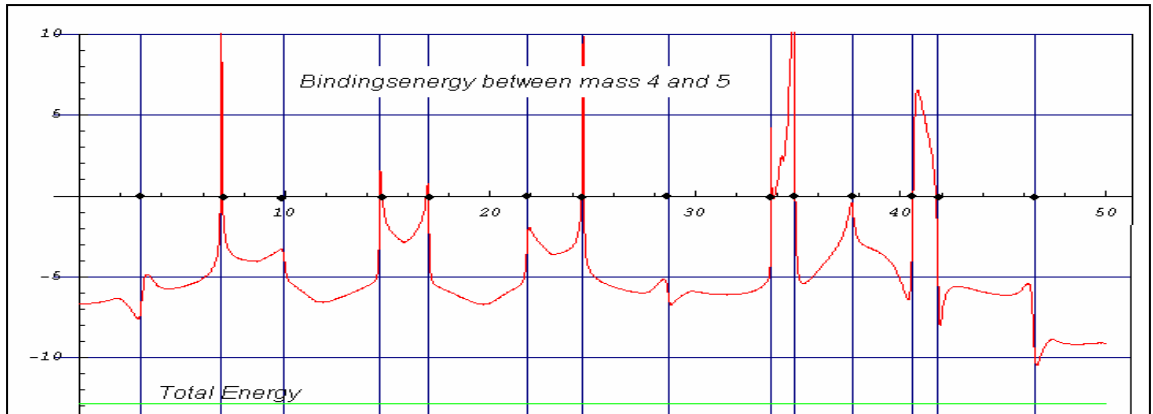

Fig. 18 -
**Bindingsenergy 1-2 (Red) ; Total Energy -14 (Green) ; ■ collinear configurations 3 in between 4,5**

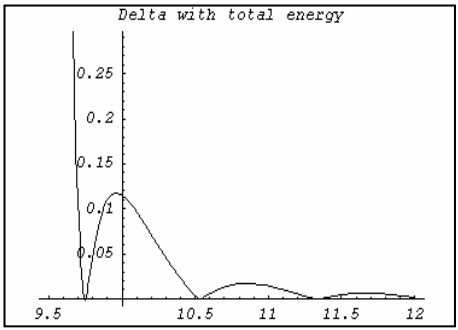
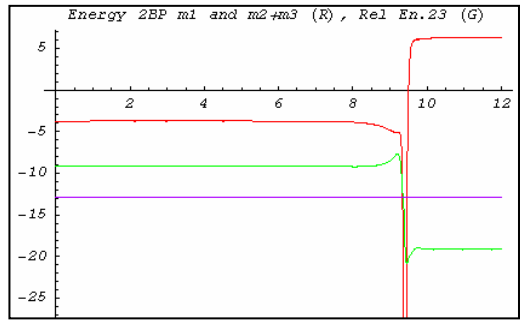

**Fig.19a-Δafter collinear configuration at t= 59.406    Fig 19b. - Crossing causing the definite split-up.**

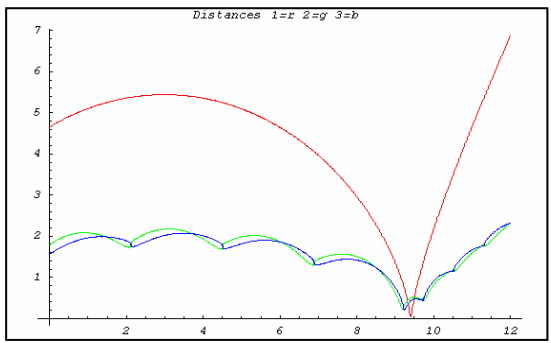

**Fig19.-Distance mass 3 from CoM-Time counted from 50.**

**Summary**
A priori, The free fall problem seems very restrictive as being planar with negative energy and zero angular momentum. However, as the zero angular momentum is a necessary condition for a triple collision the free fall problem covers this essential aspect of the 3BP.
The triple collisions occurring in central configurations can be regularized and behave as elastic collisions exactly as binary collisions in the 2 BP . The corresponding solution is periodic.  Starting from an equilateral triangle, the time to triple collision depends only on the size of the triangle and the total mass, not on the mass distribution over the particles. For the collinear central configuration, this property holds only when the outer masses are equal.
Periodic free fall solutions without any collision where the three particles come again at rest on a different triangle have been found by Standish for mass ratio's 3:4:5  in the vicinity of the pythagorean triangle with sides 3:4:5 (Burrau's problem). Standish showed also thet when the velocity of the third particle is zero at the instant of a binary collision, the solution is periodic.
Starting from a nearly equilateral triangle, illustrates  "chaos" or lack of continuity in the 3BP. The motion becomes elliptic-hyperbolic after the first close encounter. The strong interactions take place in a very short time interval about the collinear configuration.
Burrau's problem illustrates that the transition to an elliptic-hyperbolic system may occur after many close encounters. This observation combined with  accumulated numerical evidence points in the direction that, in general, a free fall motion is either periodic (with or without collisions) or becomes elliptic-hyperbolic except for the triple collisions that can not be regularized and where the motion ends

# APPENDICES

## Appendix A:  2BP summary

The standard way of solving the 2BP, is to go for the relative motion of one mass point ($m_1$) about the other ($m_2$).**[Roy]**. The result for the energy with the CoM at rest, is :

$$E = \frac{1}{2}m_1 v_1^2 + \frac{1}{2}m_2 v_2^2 - \frac{Gm_1 m_2}{|\bar{r}_2 - \bar{r}_1|}$$

in terms of  $\bar{v} = \bar{v}_1 - \bar{v}_2 \Rightarrow E = \frac{m_1 m_2}{m_1 + m_2}\left(\frac{1}{2}v^2 - \frac{Gm}{|\bar{r}|}\right) = \frac{m_1 m_2}{m_1 + m_2} E_{um}$    (A1a,b)

and  $\bar{r} = \bar{r}_1 - \bar{r}_2 \quad r = |\bar{r}| = r_1 + r_2 \quad m = m_1 + m_2$

$E_{um}$ is the energy per unit mass which has the dimensions of velocity squared. The use of this quantity is standard in the two body problem.  Here we need to distinguish clearly between $E_{um}$ and the energy of 1 or both particles.  The subscript "um" will always refers to a quantity per unit mass.

From the solution of the 2BP, we know that the relative motion of particle 1 about particle 2 is an ellipse, with particle 2 at a focus and parameters say, a (semi-major axis), p (semilatus rectum) and e (eccentricity) . The two particles describe similar ellipses with the Com, assumed at rest, in a focus. Everthing happens as if a <u>particle i is attracted by an unmovable mass $m_j^3$/m at the Com</u>  (Kepler problem) :

$$\boxed{\mu_i = G \frac{m_j^3}{m^2}}$$   (A2)

With such mass at the void origin, each particle sees the same acceleration as caused by the mutual attraction of Newton's law of gravitation.

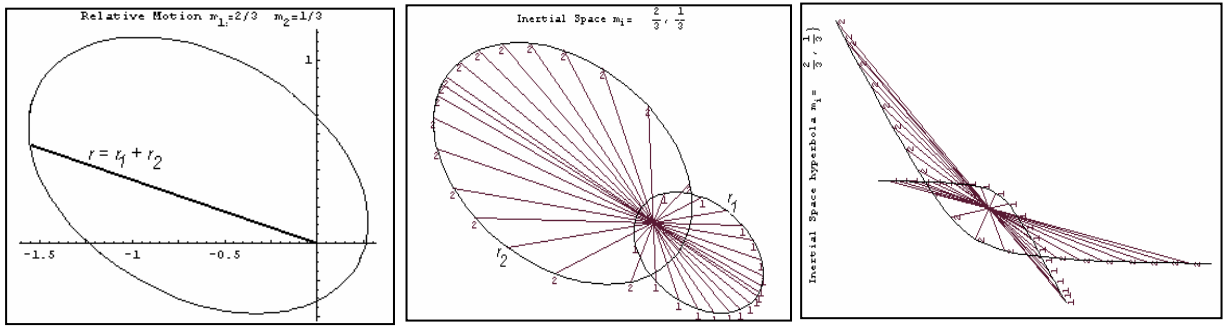

Fig.1a-c ;  2BP relative motion and motion in inertial space

So far, the Com of the 2BP was at rest. When this not the case , the potential energy is unchanged while the kinetic energy with $v_{c2}$ the velocity of the CoM,  becomes :

$$E_k = \frac{1}{2} m v_{c2}^2 + E_b$$

$$E_b = \frac{1}{2} m_1 (\bar{v}_1 - \bar{v}_c)^2 + \frac{1}{2} m(\bar{v}_2 - \bar{v}_c)^2 = \frac{1}{2}\frac{m_1 m_2}{m} v^2$$    (A3a,b)

as   $m\bar{v}_c = m_1 \bar{v}_1 + m_2 \bar{v}_2$   and   $v = \bar{v}_1 - \bar{v}_2$    as above

The energy is simply the energy of the CoM at rest augmented by the the translational energy of the CoM and **$E_b$ is called the kinetic bindingsenergy between the 2 particles.**

For the angular momentum :

$$\bar{h} = \bar{h}_1 + \bar{h}_2 = m_1 (\bar{r}_1 \times \bar{v}_1) + m_2 (\bar{r}_2 \times \bar{v}_2)$$   (4a,b)

$$= m(\bar{r}_c \times \bar{v}_c) + \frac{m_1 m_2}{m}(\bar{r} \times \bar{v})$$

The angular momentum is also simply the sum of the angular momentum of 2 masses lumped together at their the CoM and the angular momentum in their relative motion.

| Semi-major axis | | | |
|---|---|---|---|
| a = $a_1+a_2$ | $E = \dfrac{m_1 m_2}{m_1+m_2} E_{um}$ <br> $E_{um} = -\dfrac{G(m_1+m_2)}{2a}$ | $E = E_1 + E_2$ <br> $E = \dfrac{m_1 m_2}{m}\left(\dfrac{1}{2}v^2 - \dfrac{Gm}{|\bar{r}|}\right)$ <br> $E = m_1 E_{um1} + m_2 E_{um2}$ | $2a = -\dfrac{Gm}{E_{um}}$ <br> $= -\dfrac{Gm_1 m_2}{E}$ |
| $a_1 = (m_2/m)a$ | $E_1 = m_1 E_{um1}$ <br> $E_{um1} = \left(\dfrac{1}{2}v_1^2 - \dfrac{Gm_2^3}{|\bar{r}_1| m^2}\right)$ | $E_1 = m_1 \dfrac{m_2^2}{m^2}\left(\dfrac{1}{2}v^2 - \dfrac{Gm}{|\bar{r}|}\right)$ <br> $= \dfrac{m_2}{m} E$ | $2a_1 = -\dfrac{G\mu_1}{E_{um1}}$ $\mu_1 = G\dfrac{m_2^3}{m^2}$ <br> $= -\dfrac{Gm_1 m_2^3}{m^2 E_1}$ |
| $a_2 = (m_1/m)a$ | $E_2 = m_2 E_{um2}$ <br> $E_{um2} = \left(\dfrac{1}{2}v_2^2 - \dfrac{Gm_1^3}{|\bar{r}_2| m^2}\right)$ | $E_2 = m_2 \dfrac{m_1^2}{m^2}\left(\dfrac{1}{2}v^2 - \dfrac{Gm}{|\bar{r}|}\right)$ <br> $= \dfrac{m_1}{m} E$ | $2a_2 = -\dfrac{G\mu_2}{E_{um2}}$ $\mu_2 = G\dfrac{m_1^3}{m^2}$ <br> $= -\dfrac{Gm_1^3 m_2}{m^2 E_2}$ |
| Semilatus Rectum | | | |
| p = $p_1+p_2$ | $h = \dfrac{m_1 m_2}{m} h_{um}$ <br> $h_{um} = \sqrt{Gm\, p}$ | $h = h_1 + h_2$ <br> $h = m_1 h_{um1} + m_2 h_{um2}$ <br> $h = \dfrac{m_1 m_2}{m}(\bar{r}\times\bar{v})$ | $p = \dfrac{h_{um}^2}{\mu}$ $\mu = Gm$ <br> $= \dfrac{h^2}{G\dfrac{m_1 m_2}{m}} \dfrac{1}{m_1 m_2}$ |
| $P_1 = (m_2/m)p$ | $h_1 = m_1 h_{um1}$ <br> $h_{um1} = \bar{r}_1 \times \bar{v}_1$ | $h_1 = \dfrac{m_2^2}{m^2}(\bar{r}\times\bar{v})$ <br> $= \dfrac{m_2}{m} h$ | $p_1 = \dfrac{h_{um1}^2}{\mu_1}$ $\mu_1 = G\dfrac{m_2^3}{m^2}$ <br> $= \dfrac{h^2}{Gm_1^2 m_2} = \dfrac{h_1^2}{\mu_1}\dfrac{1}{m_1^2}$ |
| $P_2 = (m_1/m)p$ | $h_2 = m_2 h_{um2}$ <br> $h_{um2} = \bar{r}_2 \times \bar{v}_2$ | $h_2 = \dfrac{m_1^2}{m^2}(\bar{r}\times\bar{v})$ <br> $= \dfrac{m_1}{m} h$ | $p_2 = \dfrac{h_{um2}^2}{\mu_2}$ $\mu_2 = G\dfrac{m_1^3}{m^2}$ <br> $= \dfrac{h^2}{Gm_2^2 m_1} = \dfrac{h_2^2}{\mu_2}\dfrac{1}{m_2^2}$ |

Table I - Expressions energy and angular momentum

## Appendix B: 3BP Elliptic-Hyperbolic

Energy and Angular momentum of the 3BP with the CoM at rest :

$$E_{3BP} = \frac{1}{2}m_1 v_1^2 + \frac{1}{2}m_2 v_2^2 + \frac{1}{2}m_3 v_3^2 - \frac{Gm_1 m_2}{|\bar{r}_2-\bar{r}_1|} - \frac{Gm_1 m_3}{|\bar{r}_3-\bar{r}_1|} - \frac{Gm_3 m_2}{|\bar{r}_2-\bar{r}_3|}$$

$$\bar{H}_{3BP} = m_1(\bar{r}_1 \times \bar{v}_1) + m_2(\bar{r}_2 \times \bar{v}_2) + m_3(\bar{r}_3 \times \bar{v}_3)$$

$$m_1\bar{r}_1 + m_2\bar{r}_2 + m_3\bar{r}_3 = \mathbf{0} \quad \Rightarrow m_{12}\bar{r}_{12c} + m_3\bar{r}_3 = \mathbf{0}$$
$$m_1\bar{v}_1 + m_2\bar{v}_2 + m_3\bar{v}_3 = \mathbf{0} \quad \Rightarrow m_{12}\bar{v}_{12c} + m_3\bar{v}_3 = \mathbf{0}$$

(B1a,b)

For the elliptic-hyperbolic case we define two 2B problems : The first one has the particles $\{m_1,m_2\}$ and a moving Com. The second 2Bp has the particles $\{m_{12},m_3\}$ where $m_{12} = m_1+m_2$ located at CoM of the masses 1 and 2.

$$\frac{r_{23}}{R'} \cong \frac{m_1+m_2}{M}\sqrt{\frac{m_2+m_3}{m_2}}(1-\{quantities\ divided\ by\ R'\})$$

**Angular Momentum :**

$$\overline{H}_{3BP} = \overline{h}_1+\overline{h}_2+\overline{h}_3 = \underbrace{m_1(\overline{r}_1\times\overline{v}_1)+m_2(\overline{r}_2\times\overline{v}_2)}_{m_{12}(\overline{r}_{12c}\times\overline{v}_{12c})+h_{12r}}+m_3(\overline{r}_3\times\overline{v}_3)$$

$$= \overline{h}_{12r}+\underbrace{m_3(\overline{r}_3\times\overline{v}_3)+m_{12}(\overline{r}_{12c}\times\overline{v}_{12c})}_{\frac{m_3 m_{12}}{M}(\overline{r}_{12v}-\overline{r}_3)\times(\overline{v}_{12v}-\overline{v}_3)}$$

$$= \overline{h}_{12r}+\overline{h}_{12,3}$$

(B2a,c)

For a planar problem, $h_{12r}$ and $h_{12,3}$ have the same direction and relation 7c holds also for the moduli.

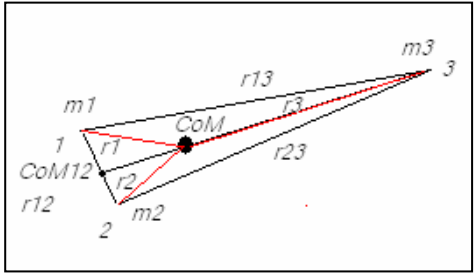

Fig 3 - Triangle after split up

Eq. B2c is exact when the CoM of the 3 particles is at rest. It does not mean that each of the two values $h_{12r}$ and $h_{12,3}$ is constant. However, they tend to a constant value when the system is split up in 2 binary systems. In particular, when the binary 1+2 and 3 is hyperbolic, the separation between particle 3 and the CoM of particles 1,2 becomes very large. The attraction on particle 3 is nearly the same as the attraction from particle 1 plus 2 placed at their CoM. Conversely, the attraction of particle 3 on the system 1,2 acts as a perturbation superposed on their mutual attraction. This perturbation has a nearly constant direction and magnitude, hence the semi-major axis or energy of the system 1,2 is not affected. The angular momentum $h_{12}$ of the particles 1,2 is w.r.t. their CoM. Everything happens as if the CoM of the system 1,2 is a particle with internal this angular momentum (spin). In the free fall case, $H_{3BP}=0$, and this "spin value" compensates the angular momentum of the hyperbolic binary system $\{3, m_1+m_2$ at their CoM$\}$

**Energy** When we group in (6a) the terms of the 2BP of particles 1,2 we have :

$$E_{3BP} = \underbrace{\frac{1}{2}m_1v_1^2+\frac{1}{2}m_2v_2^2+-\frac{Gm_1m_2}{|\overline{r}_2-\overline{r}_1|}}_{\frac{1}{2}(m_1+m_2)V_{CoM\,12}^2+E_{b12}}+\frac{1}{2}m_3v_3^2-\frac{Gm_1m_3}{|\overline{r}_3-\overline{r}_1|}-\frac{Gm_3m_2}{|\overline{r}_2-\overline{r}_3|}$$ (B3a,c)

$$= E_{b12}+\underbrace{\frac{1}{2}(m_1+m_2)V_{CoM\,12}^2+\frac{1}{2}m_3v_3^2-\frac{Gm_{12}m_3}{|\overline{r}_3-\overline{r}_{Com12}|}}_{E_{b12,3}}+\underbrace{\frac{Gm_{12}m_3}{|\overline{r}_3-\overline{r}_{com12}|}-\frac{Gm_1m_3}{|\overline{r}_3-\overline{r}_1|}-\frac{Gm_3m_2}{|\overline{r}_2-\overline{r}_3|}}_{\Delta}$$

$$\Delta = Gm_3\left[\frac{m_{12}}{|\overline{r}_3-\overline{r}_{com12}|}-\frac{m_1}{|\overline{r}_3-\overline{r}_1|}-\frac{m_2}{|\overline{r}_2-\overline{r}_3|}\right] = Gm_3\Delta'\quad and\quad R'=|\overline{r}_3-\overline{r}_{com12}|$$

R' is the distance between particle 3 and the CoM of particles 1,2. So the total energy is the sum of the bindingsenergy of the two binary systems plus $\Delta$.

$$\boxed{E_{3BP} = E_{b12}+E_{b12,3}+\Delta}$$ (B4)

When R' continues to increase R' , $|r_3-r_1|$ and $|r_3-r_2|$ become nearly equal (fig 3,) as the base 1-2 of the triangle $\{1-2-3\}$ becomes then small compared to the sides 1-3,1-2. For equal masses, R' is the median of the triangle 1-2-3. Using :

$$\overline{r}_3 = \frac{m_1+m_2}{M}(\overline{r}_3-\overline{r}_{com12})$$

we can see how R' , $|r_3-r_1|$ and $|r_3-r_2|$ become more and more equal :.

$$r_{23}^2 = \frac{m_2+m_3}{m_2}\frac{(m_1+m_2)^2}{M^2}R'^2-\left(\frac{m_1^2}{m_2m_3}r_1^2-\frac{m_2+m_3}{m_3}r_2^2\right)$$

$$r_{13}^2 = \frac{m_1+m_3}{m_1}\frac{(m_1+m_2)^2}{M^2}R'^2+\left(\frac{m_1+m_3}{m_3}r_1^2-\frac{m_2^2}{m_3m_1}r_2^2\right)$$

and a similar result for $r_{13}$. So $\lim\limits_{R' \to \infty} \dfrac{r_{23}}{R'} = \lim \dfrac{r_{13}}{R'} = \dfrac{m_1 + m_2}{M}\sqrt{\dfrac{m_2 + m_3}{m_2}}$

When the 3 masses are equal, these limits are $\sqrt{8/9}$ and of $\Delta'$ becomes : $\dfrac{m}{R'}(2 - 3/\sqrt{2}) = -\dfrac{0.12132m}{R'}$

When the quantity $\Delta'$ is small compared to the bindingsenergies, the system splits up in two binaries

$$\Delta' = \left[ \dfrac{m_1 + m_2}{R'} - \dfrac{m_1}{r_{13}} - \dfrac{m_2}{r_{23}} \right] \qquad (B5)$$

## Appendix C: On the **Rectilinear solution.**
Changing the independent variable to u : dt = ($\sqrt{a/\mu}$) x du
### a) In the EoM (Eq.7 )
$$x'' - \dfrac{x'^2}{x} + a = 0$$
It is easily verified that x-a(1+Cos u) is a solution for x(0)=2a , x'(0)=0 .
This is not obvious from this nonlinear differential equation. *Mathematica* is not able to deal with this form.
### b) In the Energy equation
$$E = \dfrac{\dot{x}^2}{2} - \dfrac{\mu}{x} \Rightarrow x' = \sqrt{x(2a-x)}$$
*Mathematica* gives the solution as : $x = 2a\dfrac{\tan^2(u/2)}{1 + \tan^2(u/2)}$ which can be simplified by hand to
2a Sin$^2$(u/2)  and finally a(1+Cos u).

## Appendix D : Transformation of the isosceles case
To regularize the BC's => r= y in the relation between du and dt is indicated , while ;for the TC's r = ρ seems appropriate. The coupling between the equations 13a-b makes that neither of these choices is useful. Another interesting length is the radius of gyration $r_g$: which is equivalent with the out of plane inertia $J_z$ :

$$M\, r_g^2 = J_z = \sum_{i=1,3} m_i r_i^2 = \sum_{i=1,3, \neq j} \dfrac{m_i m_j}{M} r_{ji}^2 \qquad (19)$$

In our case : $\quad r_g = \sqrt{\dfrac{2}{3}}\sqrt{\dfrac{x^2}{3} + y^2} = \dfrac{\sqrt{2}}{3}\sqrt{\rho^2 + 2y^2} \qquad (20)$

With the radius of gyration $r_g$ we can associate a direction defined by an angle $\theta_c$ :

$$\cos\theta_c = \dfrac{\sqrt{2}}{\sqrt{3}}\dfrac{y}{r_g} \quad \sin\theta_c = \dfrac{\sqrt{2}}{3}\dfrac{x}{r_g} \quad => \tan\theta_c = \dfrac{1}{\sqrt{3}}\dfrac{x}{y} \quad \text{and}\; \tan\theta_b = \dfrac{x}{y} \qquad (21)$$

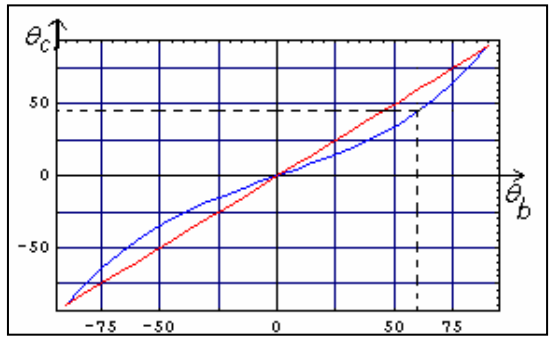

where $\theta_b$ is the base angle of the isosceles triangle :
{ $-\pi/2 \le \theta_b \le \pi/2$}. The relation between $\theta_v$ and $\theta_c$ is shown in Fig.7. Note that $\theta_b$=60° (equilateral), corresponds to $\theta_c$=45°. The angle $\theta_o$ from the CoM to the base particle (y > 0) is $\theta_o = \theta_b + \pi/2$ { $0 \le \theta_o \le \pi$}. A configuration in {y, x} can now be visualised by 3 particles at $r_g$, : 1 on the x-axis and the 2 others making an angle ± $\theta_c$ with the axis (Fig.8)

**Fig.7 Relation θ_v and θ_c**

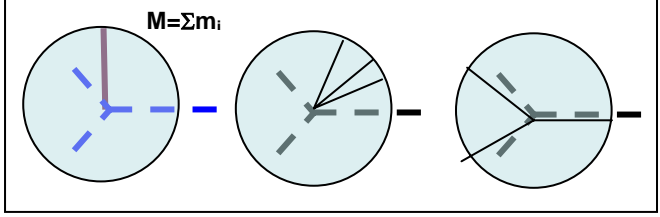

**Fig.8 - Radius of gyration ; θ_c**

With the point transformation :
$$x = \frac{3}{\sqrt{2}} r_g \sin\theta_c \quad y = \frac{\sqrt{3}}{\sqrt{2}} r_g \cos\theta_c \tag{20}$$

the energy per unit mass $E_u = E/m$ eq. 15c becomes :

$$E_u = T_u + V_u \quad \text{with} \quad T_u = \frac{3}{2}\left(\dot{r}_g^2 + r_g^2 \dot{\theta}_c^2\right) \quad V_u = -\frac{\mu}{\sqrt{6}\, r_g}\left(\frac{1}{\cos\theta_c} + \frac{4}{\sqrt{1+2\sin^2\theta_c}}\right) \tag{21}$$

$V_u$ is now written as a product of a function of $r_g$ and $F(\theta_c)$ depending only on $\theta_c$

$$V_u = -\frac{\mu}{\sqrt{6}\, r_g} F(\theta_c) \Rightarrow F(\theta_c) = \left(\frac{1}{\cos\theta_c} + \frac{4}{\sqrt{1+2\sin^2\theta_c}}\right) \tag{22}$$

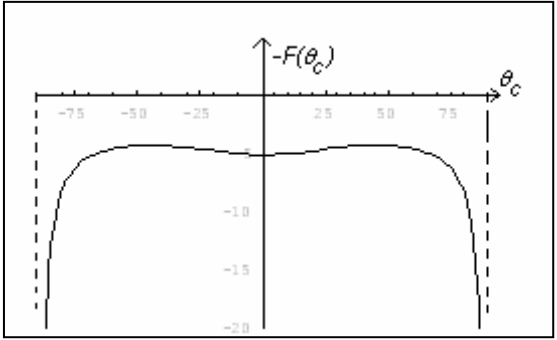

The function $F(\theta_c)$ is shown in Fig.9 and plays an important role in the analysis of the singularities. When revolving about the $q_c$ -axis, it generates a cylindrical surface known as the triple collision manifold.

**Fig.9 - Radius of gyration ; θ_c**

### Appendix E : Neighbouring solutions Equilateral triangle

Starting from a nearly equilateral triangle it is possible to write out the perturbation equations as the exact solution to a fall from an exact equilateral triangle is known [eq.10 -12] . The starting point are the perturbation equations to the straight line solution of the 2BP in the form given in eqn.10 :

$$\ddot{x}_1 = -\frac{2\mu}{r_r^3} x_1 \quad r_r^2 = x_r^2 + y_r^2 = x_r^2 + 0$$
$$\ddot{y}_1 = -\frac{\mu}{r_r^3} y_1 \quad . \tag{E1-a,b}$$

where $r_r$ is the relative distance between the point in the straight line solution.

$x_1\ y_1$ are the relative first order deviations from the straight line at the same time.

This equation must be converted to the deviations of one point to the CoM with the appropriate μ value as given by eq.12a-c). For the extension to the 3BP, all deviations must be expressed in the same frame. Starting from the general 3BP equations the linearized deviations of one point will be coupled to the deviations of the two other points. The linearized system has periodic coefficients.

The result for the 2BP after changing to the independent variable u, gives the following two decoupled equations :

$$x_1'' + \tan\frac{u}{2} x_1' - \frac{x}{\cos^2(u/2)} = 0 \qquad \text{(E2-a,b)}$$
$$y_1'' + \tan\frac{u}{2} y_1' + \frac{y}{2\cos^2(u/2)} = 0$$

The fundamental solutions are :

for the y- equation : $e_{1y}(u) = \sin u \qquad e_{1y}(0) = 0 \quad e_{1y}'(0) = 1$

$e_{2y}(u) = (1+\cos u)/2 \quad e_{2y}(0) = 1 \quad e_{2y}'(0) = 0$

for the x -equation : $e_{1x}(u) = 2\tan u/2 \qquad\qquad\qquad e_{1x}(0) = 0 \quad e_{1x}'(0) = 1$

$e_{2x} = \dfrac{1}{5\cdot 2} F_1(1,2,7/2;\cos^2 u/2)\cos^4 u/2 \quad e_{2x}(0) = 1 \quad e_{2x}'(0) = ??$